\documentclass[final,5p,times,8pt]{elsarticle}
\usepackage[centertags]{amsmath}
\usepackage{amsfonts,amssymb,amsthm,epstopdf,newlfont,graphicx}
\usepackage{booktabs,bm,enumitem,setspace,empheq,cancel}
\usepackage{algorithm}
\usepackage[short,nocomma]{optidef}
\usepackage{algpseudocode}
\usepackage{caption,subcaption}
\usepackage{threeparttable,multirow}
\usepackage{amsmath,accents}
\usepackage[T1]{fontenc}
\usepackage[american]{babel}
\usepackage{natbib}
\biboptions{comma,square,sort&compress}
\usepackage{appendix}
\usepackage[colorlinks=true]{hyperref}

\newcommand\hcancel[2][black]{\setbox0=\hbox{$#2$}%
\rlap{\raisebox{.25\ht0}{\textcolor{#1}{\rule{0.7\wd0}{0.75pt}}}}#2} 
\newcommand\hcancelt[2][black]{\setbox0=\hbox{$#2$}%
\rlap{\raisebox{.25\ht0}{\textcolor{#1}{\hspace{0.3mm}\rule{0.7\wd0}{0.75pt}}}}#2} 
\newtheorem{thm}{Theorem}[section]

\newtheorem{rem}{Remark}[section]
\theoremstyle{definition}

\numberwithin{algorithm}{section}
\numberwithin{equation}{section}
\renewcommand{\theequation}{\thesection.\arabic{equation}}
\def\simgt{\,\hbox{\lower0.6ex\hbox{$>$}\llap{\raise0.3ex\hbox{$\sim$}}}\,}
\def\simlt{\,\hbox{\lower0.6ex\hbox{$<$}\llap{\raise0.3ex\hbox{$\sim$}}}\,}
\def\simgteq{\,\hbox{\lower0.6ex\hbox{$\ge$}\llap{\raise0.6ex\hbox{$\sim$}}}\,}
\def\simlteq{\,\hbox{\lower0.6ex\hbox{$\le$}\llap{\raise0.6ex\hbox{$\sim$}}}\,}
\def\applteq{\,\hbox{\lower0.6ex\hbox{$\le$}\llap{\raise0.8ex\hbox{$\approx$}}}\,}
\def\applt{\,\hbox{\lower0.6ex\hbox{$<$}\llap{\raise0.5ex\hbox{$\approx$}}}\,}
\DeclareMathAlphabet\mathbfcal{OMS}{cmsy}{b}{n}

\DeclareMathOperator*{\indmax}{indmax}
\DeclareMathOperator*{\indmin}{indmin}

\DeclareMathOperator{\dist}{dist}

\makeatletter
\def\user@resume{resume}
\def\user@intermezzo{intermezzo}
\newcounter{previousequation}
\newcounter{lastsubequation}
\newcounter{savedparentequation}
\renewenvironment{subequations}[1][]{%
      \def\user@decides{#1}%
      \setcounter{previousequation}{\value{equation}}%
      \ifx\user@decides\user@resume 
           \setcounter{equation}{\value{savedparentequation}}%
      \else  
      \ifx\user@decides\user@intermezzo
           \refstepcounter{equation}%
      \else
           \setcounter{lastsubequation}{0}%
           \refstepcounter{equation}%
      \fi\fi
      \protected@edef\theHparentequation{%
          \@ifundefined {theHequation}\theequation \theHequation}%
      \protected@edef\theparentequation{\theequation}%
      \setcounter{parentequation}{\value{equation}}%
      \ifx\user@decides\user@resume 
           \setcounter{equation}{\value{lastsubequation}}%
         \else
           \setcounter{equation}{0}%
      \fi
      \def\theequation  {\theparentequation  \alph{equation}}%
      \def\theHequation {\theHparentequation \alph{equation}}%
      \ignorespaces
}{%
  \ifx\user@decides\user@resume
       \setcounter{lastsubequation}{\value{equation}}%
       \setcounter{equation}{\value{previousequation}}%
  \else
  \ifx\user@decides\user@intermezzo
       \setcounter{equation}{\value{parentequation}}%
  \else
       \setcounter{lastsubequation}{\value{equation}}%
       \setcounter{savedparentequation}{\value{parentequation}}%
       \setcounter{equation}{\value{parentequation}}%
  \fi\fi
  \ignorespacesafterend
}
\makeatother
\newcommand{\C}[1]{\mathcal{#1}}

\newcommand{\F}[1]{\mathbf{#1}}

\newcommand{\MB}[1]{\mathbb{#1}}
\newcommand{\MBS}{\MB{S}}
\newcommand{\MBR}{\mathbb{R}}
\newcommand{\MBE}{\mathbb{E}}
\newcommand{\MBRP}{\MBR^+}

\newcommand{\MBZ}{\mathbb{Z}}
\newcommand{\MBZP}{\MBZ^+}
\newcommand{\MBZzer}{\MBZ_0}
\newcommand{\MBZzerP}{\MBZzer^+}
\newcommand{\MBZe}{\MBZ_e}
\newcommand{\MBZeP}{\MBZe^+}

\newcommand{\MBT}{\mathbb{T}}
\newcommand{\MBJ}{\mathbb{J}}
\newcommand{\MBN}{\mathbb{N}}

\newcommand{\MBC}{\mathfrak{C}}
\newcommand{\MBK}{\mathfrak{K}}
\newcommand{\FThe}{\F{\Theta}}
\newcommand{\FThetpi}{\FThe^{2\pi}}

\newcommand{\cancbra}[1]{\hcancel{[}#1\hcancelt{]}}
\newcommand{\sumd}{\sideset{}{'}}

\newcommand{\bmx}{\bm{x}}

\newcommand{\bmt}{\bm{t}}
\newcommand{\bmtau}{\bm{\tau}}
\newcommand{\bmT}{\bm{T}}

\newcommand{\bmz}{\bm{z}}
\newcommand{\bmD}{\bm{D}}
\newcommand{\bmV}{\bm{V}}
\newcommand{\bmgam}{\bm{\gamma}}
\newcommand{\bmalp}{\bm{\alpha}}
\newcommand{\bmxi}{\bm{\xi}}

\newcommand{\tilx}{\tilde{x}}
\newcommand{\tilxi}{\tilde{\xi}}

\newcommand{\tilz}{\tilde{z}}
\newcommand{\tilD}{\tilde{D}}
\newcommand{\tilT}{\tilde{T}}
\newcommand{\tilV}{\tilde{V}}
\newcommand{\tilgam}{\tilde{\gamma}}
\newcommand{\tilalp}{\tilde{\alpha}}
\newcommand{\tilbmx}{\tilde{\bmx}}

\newcommand{\tilbmz}{\tilde{\bmz}}
\newcommand{\tilbmD}{\tilde{\bmD}}
\newcommand{\tilbmT}{\tilde{\bmT}}
\newcommand{\tilbmV}{\tilde{\bmV}}
\newcommand{\tilbmgam}{\tilde{\bmgam}}
\newcommand{\tilbmalp}{\tilde{\bmalp}}
\newcommand{\tilbmxi}{\tilde{\bmxi}}
\newcommand{\barJ}{\bar{J}}
\newcommand{\bmy}{\bm{y}}

\newcommand{\bmzer}{\bm{\mathit{0}}}
\newcommand{\bmone}{\bm{\mathit{1}}}

\newcommand{\bmX}{\bm{X}}

\newcommand{\Nin}{N_{\text{in}}}
\newcommand{\Ninc}{N_{\text{inc}}}

\newcommand{\hxi}{\hat{\xi}}

\newcommand{\FOmega}{\F{\Omega}}
\newcommand{\rme}{\mathrm{e}_0}
\newcommand{\foralla}{\,\forall_{\mkern-6mu a}\,}
\newcommand{\foralle}{\,\forall_{\mkern-6mu e}\,}
\newcommand{\foralls}{\,\forall_{\mkern-6mu s}\,}
\newcommand{\foralll}{\,\forall_{\mkern-4mu l}\,}

\usepackage{mathtools}
\mathtoolsset{showonlyrefs}
\makeatletter
\def\BState{\State\hskip-\ALG@thistlm}
\makeatother
\MakeRobust{\eqref}
\errorcontextlines\maxdimen

\makeatletter
    \newcommand*{\algrule}[1][\algorithmicindent]{\makebox[#1][l]{\hspace*{.5em}\thealgruleextra\vrule height \thealgruleheight depth \thealgruledepth}}%
\newcommand*{\thealgruleextra}{}
\newcommand*{\thealgruleheight}{.75\baselineskip}
\newcommand*{\thealgruledepth}{.25\baselineskip}

\newcount\ALG@printindent@tempcnta
\def\ALG@printindent{%
    \ifnum \theALG@nested>0
        \ifx\ALG@text\ALG@x@notext
        \else
            \unskip
            \addvspace{-1pt}
            \ALG@printindent@tempcnta=1
            \loop
                \algrule[\csname ALG@ind@\the\ALG@printindent@tempcnta\endcsname]%
                \advance \ALG@printindent@tempcnta 1
            \ifnum \ALG@printindent@tempcnta<\numexpr\theALG@nested+1\relax
            \repeat
        \fi
    \fi
    }%
\usepackage{etoolbox}
\patchcmd{\ALG@doentity}{\noindent\hskip\ALG@tlm}{\ALG@printindent}{}{\errmessage{failed to patch}}
\makeatother

\newbox\statebox
\newcommand{\myState}[1]{%
    \setbox\statebox=\vbox{#1}%
    \edef\thealgruleheight{\dimexpr \the\ht\statebox+1pt\relax}%
    \edef\thealgruledepth{\dimexpr \the\dp\statebox+1pt\relax}%
    \ifdim\thealgruleheight<.75\baselineskip
        \def\thealgruleheight{\dimexpr .75\baselineskip+1pt\relax}%
    \fi
    \ifdim\thealgruledepth<.25\baselineskip
        \def\thealgruledepth{\dimexpr .25\baselineskip+1pt\relax}%
    \fi
    \State #1%
    \def\thealgruleheight{\dimexpr .75\baselineskip+1pt\relax}%
    \def\thealgruledepth{\dimexpr .25\baselineskip+1pt\relax}%
}

\begin{document}
\begin{frontmatter}
\title{Optimal Periodic Control of Unmanned Aerial Vehicles Based on Fourier Integral Pseudospectral and Edge-Detection Methods}
\author[Assiut]{Kareem T. Elgindy\corref{cor1}}
\ead{kareem.elgindy@(aun.edu.eg;gmail.com)}
\address[Assiut]{Mathematics Department, Faculty of Science, Assiut University, Assiut 71516, Egypt}
\cortext[cor1]{Corresponding author}
\begin{abstract}
This study describes the development of a novel numerical optimization framework to maximize the endurance of unmanned aerial vehicles (UAVs). We address the problem of numerically determining the optimal thrust and cruise angle of attack in a two-dimensional space for a UAV under certain initial, periodic, and bound constraints. The time horizon of the free final time optimal control problem (OCP) is first normalized, and the normalized OCP in integral form is discretized in physical space into a nonlinear programming problem (NLP) using Fourier collocation and quadrature based on equispaced points. Great attention in this work is placed on the accurate detection of jump discontinuities and resolving the thrust history effectively directly from the Fourier pseudospectral (FPS) data through a novel edge-detection method without any smoothing techniques. The numerical results demonstrate that the proposed method is simple, stable, and easy to implement.
\end{abstract}

\end{frontmatter}

\section{Introduction}
\label{Int}
Unmanned aerial vehicles (UAVs) have gained significant importance in recent decades in military, commercial, and civilian applications, such as reconnaissance, target engagement, combat support, traffic and border control, aerial surveillance, law enforcement, remote parcel delivery, hidden and hazardous area exploration, onset detection of subpavement voids, maintenance and repair of aircraft, crop spraying, wildfire fighting, search and rescue operations, climate change monitoring, etc. \cite{mobariz2015long,kim2015situation,wang2022evolutionary,mechan2023unmanned,
mahmud2023unmanned,kulkarni2023deep}. Part of the success of these vehicles is commonly due to their (i) low prices, (ii) low cost utility, (iii) low maintenance, (iv) faster deployment, (v) easier, efficient, and practical operation compared with conventional manned aircraft and road delivery vehicles, (vi) large capacity to take closer footage without compromising the quality of both photos and video, (vii) reduction of pilot fatality and injury rates mid-flights, (viii) ability to land or take off without having to use runways, etc. \cite{cotter2019application,meng2021space,wang2022steady,
gunaratne2022unmanned,borowik2022mutable}. 

To cover long ranges, UAV industries are in dire need of extended flight endurance. Increasing the amount of fuel UAVs can carry will increase their total weight and reduce their flight endurance. Fortunately, several other approaches exist to extend the endurance and allow for long continuous flights. One approach is to use certain fueled systems, such as liquid hydrogen- or hydrocarbon-fueled systems \cite{Sweetman2006,de2020design}, although the latter systems generally suffer from fundamental problems pertaining to noise, efficiency, and reliability \cite{khofiyah2019technical}. Another approach is to use solar cells to power UAVs by harvesting energy from the outside and providing additional power that is dependent on the sun, weight, wing surface area, and efficiency \cite{scheiman2016path,dwivedi2018maraal,mateja2023energy}. Other research works were stirred on developing more efficient propellers as a viable approach \cite{nguyen2015possibility,yonezawa2016propeller,vijayanandh2019research,
Wisniewski2022,yang2023high}. One of the cheapest yet most effective approaches to considerably boost the endurance and performance of UAVs is to minimize the amount of fuel consumption required for flight operation through optimal trajectory planning \cite{dobrokhodov2020energy,
dobrokhodov2020achievable,wenjun2022energy,xi2022energy}. This is commonly achieved by optimally regulating the thrust produced by the UAV propellers, attack angle, and bank angle to realize autonomous flight based on energy maximization. Closely related works in this direction investigated the possibility of finding optimal ``periodic'' solutions to enhance endurance rather than steady-state optimal solutions; \cite{sachs2009unlimited,hosseini2013optimal,wenkai2017optimal,
wang2019periodically,ogunbodede2019periodic,ogunbodede2019endurance}. This interesting subject is part of the optimal periodic control theory, which was originally motivated by problems from chemical engineering, as some studies have found that cycling a chemical reactor can increase the average output compared to steady-state operation \cite{colonius2006optimal}. Later, optimal periodic control theory found its way in performance optimization of satellites, aircraft flights, ships and passenger cars during cruising, diesel engines, bioreactors, drug delivery, etc. \cite{gilbert1976vehicle,higuchi2010optimal,ghanaatpishe2017structure,
sivertsson2017optimal,shen2018fuel,Elgindy2023a}.

In this study, we are interested in the periodic energy-optimal path planning of UAVs. We present the Fourier integral pseudospectral method integrated with an edge-detection technique (FIPS-ED method): A novel Fourier integral pseudospectral method (FIPS)-based direct optimization method integrated with a robust edge-detection technique to determine the optimal thrust and flight angle of attack required to maximize the UAV endurance by minimizing the rate of fuel consumption per unit time. The proposed method combines a time-scaling strategy with integral reformulations to convert the problem into a normalized optimal control problem (OCP) in integral form. Fourier collocation and quadrature induced by the accurate and efficient Fourier integration matrix (FIM) are then applied to discretize the problem successfully into a constrained nonlinear programming problem (NLP), which can be treated using standard optimization software. A novel edge-detection technique is also introduced to accurately locate the discontinuity points of the bang-bang thrust and effectively reconstruct it directly from the Fourier pseudospectral (FPS) data; the proposed technique, which we call the FPSED method, is a robust development of an earlier version which appeared recently in \cite{Elgindy2023a}. We show further that smoothing techniques combined with Fourier-based methods of low mode number or mesh densities do not generally simulate 2D-dimensional UAV flights accurately, because they can smear the discontinuities of bang-bang thrusts, and the errors in the designed thrust may cause deviations in the optimal state variables of the flight. Instead, an adequate number of Fourier modes or mesh points are necessary for numerical optimization methods to properly derive the optimal thrust policy while maintaining the validity of the UAV model. We assume in this work that the flight simulation is modeled by the widely accepted 2D point-mass dynamic model subject to certain initial, periodic, and bound constraints on the state and control variables, and that the final time of the flight duration is free. Moreover, the extreme thrust values are given, but the thrust optimal switching times are not known a priori. For a comprehensive survey on the excellent virtues of integral reformulations and FIMs used in our work, we refer to \cite{du2016well,elgindy2018highb,elgindy2019high,elgindy2020distributed,
Elgindy2023a,Elgindy2023b} and the references therein. Clear expositions of the IPS methods can also be found in \cite{elgindy2013solving,elgindy2013fast,tang2016new,Elgindy20171,
dahy2021high,ElgindyHareth2023a} and the references therein.

The rest of the paper is organized as follows: Section \ref{sec:PN} introduces some preliminary notations used in the paper. Section \ref{sec:PS1} describes the 2D path planning problem under study. In Section \ref{sec:FIPSQ1}, we review and derive some useful FPS interpolation and quadrature formulas  pertinent to the forthcoming development. Section \ref{sec:FPIMIRF1} presents the proposed FIPS-ED method. A prescription of the FPSED method is discussed in Section \ref{subsec:TFEDM1}. In Section \ref{sec:CRFSPF1}, we state the main errors and convergence results associated with the employed numerical tools. Simulation results are shown in Section \ref{sec:NS1}. Finally, we conclude the paper with some remarks in Section \ref{sec:Conc}.

\section{Preliminary Notations}
\label{sec:PN}
The following notations are used throughout this paper to abridge and simplify the mathematical formulas. Many of these notations appeared earlier in \cite{Elgindy2023a,Elgindy2023b}; however, for convenience and to keep the paper self-explanatory, we summarize them below together with the new notations.

\noindent\textbf{Logical Symbols.} The  symbols $\forall, \foralla, \foralle, \foralls$, and $\foralll$ stand for the phrases ``for all,'' ``for any,'' ``for each,'' ``for some,'' and ``for a relatively large'' in respective order. $:=$ means ``is replaced by'' or ``updated with'' and we use it here when we need to update the value of a certain variable or set. For example, $x := x + 1$ means $x$ is replaced by its old value plus one.\\[0.5em]
\textbf{List and Set Notations.} $\MBC$ denotes the set of all complex-valued functions. Moreover, $\MBR, \MBZ, \MBZP, \MBZzerP$, and $\MBZeP$ denote the sets of real numbers, integers, positive integers, non-negative integers, and positive even integers, respectively. The notations $i:j:k$ or $i(j)k$ indicate a list of numbers from $i$ to $k$ with increment $j$ between numbers, unless the increment equals one where we use the simplified notation $i:k$. For example, $0:0.5:2$ simply means the list of numbers $0, 0.5, 1, 1.5$, and $2$, while $0:2$ means $0, 1$, and $2$. The list of symbols $y_1, y_2, \ldots, y_n$ is denoted by $\left. y_i \right|_{i=1:n}$ or simply $y_{1:n}$, and their set is represented by $\{y_{1:n}\}\,\foralla n \in \MBZP$. We define $\MBJ_n = \{0:n-1\}, \MBJ_n^+ = \MBJ_n \cup \{n\}, \MBN_n = \{1:n\}\,\foralla n \in \MBZP$, and $\MBK_N = \{-N/2:N/2\}\,\foralla N \in \MBZeP$. Also, $\MBS_n^{\C{T}} = \left\{t_{0:n-1}\right\}$ is the set of $n$ equally-spaced points such that $t_j = \C{T} j/n\, \forall j \in \MBJ_n$; in this case, we write $t_n = \C{T}$. The specific interval $[0, \C{T}]$ is denoted by $\FOmega_{\C{T}}\,\forall \C{T} > 0$. For example, $[0, t_{n}]$ is denoted by ${\FOmega_{t_{n}}}$. The notation ${\FOmega_{t_{0:N-1}}}$ stands for the list of intervals ${\FOmega_{t_{0}}}, {\FOmega_{t_{1}}}, \ldots, {\FOmega_{t_{N-1}}}$; moreover, the notations ${}_a\FOmega_b$ and $\left. {{}_{t_i}\FOmega_{t_{i+1}}} \right|_{i=0:n-1}$ denote $[a, b]$ and the list of intervals $[t_0, t_1],  \ldots, [t_{n-1}, t_n]$ in respective order. $\bm{\beta} = [-\beta, \beta]\,\forall \beta > 0$, and ${\F{C}_{\C{T},\beta} } = \left\{ {x + iy:x \in {\FOmega_{\C{T}}},y \in \bm{\beta}} \right\}\;\forall \beta  > 0$. $\dist(\F{A},\F{B}) = \inf\{|x-y|: x \in \F{A}, y \in \F{B}\}\,\foralla$ nonempty sets $\F{A}, \F{B} \subseteq \MBR$. Finally, ${f^\to}(\F{A})$ denotes the image set of a function $f$ defined on a set $\F{A}.$\\[0.5em] 
\textbf{Function Notations.} $\delta_{n,m}$ is the usual Kronecker delta function of variables $n$ and $m$. For convenience, we shall denote $g(t_{n})$ by $g_n \foralla g \in \MBC, n \in \MBZ, t_n \in \MBR$, unless stated otherwise.\\[0.5em]
\textbf{Integral Notations.} We denote $\int_0^{b} {h(t)\,dt}$ by $\C{I}_{b}^{(t)}h \foralla$ integrable $h \in \MBC, b \in \MBR$. If the integrand function $h$ is to be evaluated at any other expression of $t$, say $u(t)$, we express $\int_0^{b} {h(u(t))\,dt}$ with a stroke through the square brackets as $\C{I}_{b}^{(t)}h\cancbra{u(t)}$.\\[0.5em] 
\textbf{Space and Norm Notations.} $\MBT_{\C{T}}$ is the space of $\C{T}$-periodic, univariate functions $\foralla \C{T} \in \MBRP$. $C^k(\FOmega)$ is the space of $k$ times continuously differentiable functions on ${\FOmega}\,\forall k \in \MBZzerP$.  $L^p({\FOmega_{\C{T}}})$ is the Banach space of measurable functions $u$ defined on ${\FOmega_{\C{T}}}$ such that ${\left\| u \right\|_{{L^p}}} = {\left( {{\C{I}_{\FOmega_{\C{T}}}}{{\left| u \right|}^p}} \right)^{1/p}} < \infty\,\forall p \ge 1$. The space
\[\displaystyle{{H^s}({\FOmega_{\C{T}}}) = \left\{ {u \in {L_{loc}}({\FOmega_{\C{T}}}),\;{D^\alpha }u \in {L^2}({\FOmega_{\C{T}}}),\left| \alpha  \right| \le s\;} \right\}}\quad \forall s \in \MBZzerP,\]
is the inner product space with the inner product 
\[{(u,v)_s} = \sum\nolimits_{\left| \alpha  \right| \le s} {\C{I}_{{\FOmega_{\C{T}}}} ^{(t)}\left( {{D^\alpha }u\,{D^\alpha }v} \right)},\] where ${{L_{loc}}({\FOmega_{\C{T}}} )}$ is the space of locally integrable functions on ${\FOmega_{\C{T}}}$ and ${{D^\alpha }u}$ denotes any derivative of $u$ with multi-index $\alpha$. Moreover,
\[\C{H}_{\C{T}}^s = \left\{ {u \in {H^s}({\FOmega_{\C{T}}}),\;{u^{(s)}} \in {BV},\;{u^{(0:s - 1)}}(0) = {u^{(0:s - 1)}}(\C{T})} \right\},\]
where $u^{(0:s-1)}$ denotes the column vector of derivatives $[u, u'$, $\ldots, u^{(s-1)}]^{\top}$, and $\displaystyle{{BV} = \left\{ {u \in {L^1}({\FOmega_{\C{T}}}):{{\left\| u \right\|}_{BV}} < \infty } \right\}}$ with the norm 
\[{{\left\| u \right\|}_{BV}} = \sup \left\{ {\C{I}_{\C{T}}^{(x)}(u\phi '),\;\phi  \in \C{D}({\FOmega_{\C{T}}}),\;{{\left\| \phi  \right\|}_{{L^\infty }}} \le 1} \right\}\]
such that $\C{D}({\FOmega_{\C{T}}}) = \left\{u \in {C^\infty }({\FOmega_{\C{T}}}):{\text{supp}}(u)\text{ is a compact subset }\right.$ of $\left. {\FOmega_{\C{T}}} \right\}$. The space $\C{B}_{\C{T}} = \{u \in \C{H}_{\C{T}}^0: {u^\to}\left(\FOmega_{\C{T}}\right) = \MBE_u, u(0) = u(\C{T})\}$ is the space of $\C{T}$-periodic bang-bang\footnote{A bang-bang function is a two-state, piecewise constant function that switches abruptly between two states.} functions, where $\MBE_u = \{u_1, u_2 \in \MBR: u_1 \ne u_2\}$. For analytic functions, we define the space\\
\scalebox{0.93}{$\C{A}_{\C{T},\beta} = \{u \in \C{H}_{\C{T}}^{\infty}: u \text{ is analytic in some open set containing }\F{C}_{\C{T},\beta}\},$}\\
with the norm ${\left\| u \right\|_{{\C{A}_{\C{T},\beta}}}} = {\left\| u \right\|_{{L^\infty }({{\F{C}}_{\C{T},\beta}})}}$. For convenience of writing, we shall denote ${\left\|  \cdot  \right\|_{{L^2}({\FOmega_{\C{T}}})}}$ by $\left\|  \cdot  \right\|$, and call a function $u \in \C{A}_{\C{T},\beta}$ \textit{``a $\beta$-analytic function''} if $u$ is analytic on ${\F{C}_{\C{T},\infty}}$ and $\displaystyle{{\lim _{\beta  \to \infty }}\frac{{{{\left\| u \right\|}_{{\C{A}_{\C{T},\beta} }}}}}{{{e^{{\omega _\beta }}}}} = 0}$, where ${\omega _{\alpha}} = 2\pi \alpha/\C{T}\,\forall \alpha \in \MBR$.\\[0.5em]
\textbf{Vector Notations.} We shall use the shorthand notation $\bmt_N$ and $g_{0:N-1}$ to stand for the column vectors $[t_{0}, t_{1}, \ldots$, $t_{N-1}]^{\top}$ and $[g_0, g_1, \ldots, g_{N-1}]^{\top}\,\forall N \in \MBZP$ in respective order. $\bmt_{n:m}$ denotes the subvector of $\bmt_N$ containing all elements from $t_n$ to $t_m$ in ascending order $\foralla n, m \in \MBJ_N: n < m$. The notation ${(\bm{y})_{(N)}}$ stands for $\underbrace {\bm{y} \odot \bm{y} \odot \ldots \odot \bm{y}}_{N - {\text{times}}}\; \foralla$ vector $\bmy$, where $\odot$ denotes the Hadamard (entrywise) product. In general, $\foralla h \in \MBC$ and vector $\bmy$ whose $i$th-element is $y_i \in \MBR$, the notation $h(\bmy)$ stands for a vector of the same size and structure of $\bmy$ such that $h(y_i)$ is the $i$th element of $h(\bmy)$. Furthermore, we adopt the notation $\C{I}_{{\bmt_{N}}}^{(t)}h$ to denote the $N$th-dimensional column vector $\left[ {\C{I}_{{t_{0}}}^{(t)}h,\C{I}_{{t_{1}}}^{(t)}h, \ldots ,\C{I}_{{t_{N - 1}}}^{(t)}h} \right]^{\top}$.\\[0.5em] 
\textbf{Matrix Notations.} $\F{O}_n, \F{1}_n$, and $\F{I}_n$ stand for the zero, all ones, and the identity matrices of size $n$. $\F{C}_{n,m}$ indicates that $\F{C}$ is a rectangular matrix of size $n \times m$; moreover, $\F{C}_n$ denotes a row vector whose elements are the $n$th-row elements of $\F{C}$, except when $\F{C}_n = \F{O}_n, \F{1}_n$, or $\F{I}_n$, where it denotes the size of the matrix. For convenience, a vector is represented in print by a bold italicized symbol while a two-dimensional matrix is represented by a bold symbol, except for a row vector whose elements form a certain row of a matrix where we represent it in bold symbol as stated earlier. For example, $\bmone_n$ and $\bmzer_n$ denote the $n$-dimensional all ones- and zeros- column vectors, while $\F{1}_n$ and $\F{O}_n$ denote the all ones- and zeros- matrices of size $n$, respectively. Finally, the notation $[.;.]$ denotes the usual vertical concatenation.

\section{Problem Statement}
\label{sec:PS1}
The OCP involves finding the optimal periodic attack angle $\alpha$ and thrust $T$ and their corresponding positions of the aircraft centre of gravity in the Flat-Earth reference frame $x$ and $z$, the flight-path angle $\gamma$, and the aircraft speed $V$ in the time interval $\FOmega_{T_f}$ that maximizes the UAV endurance by minimizing the rate of fuel consumption per unit time
\begin{subequations}
\begin{equation}\label{eq:eq1}
J = \frac{\sigma}{T_f} \C{I}_{T_f}^{(t)} {T}
\end{equation}
subject to the 2D point-mass UAV dynamic model 
\begin{alignat}{3}
\dot x &= V \cos \gamma,\quad &&\dot z &&= V \sin \gamma,\label{eq:sysDE1}\\
\dot \gamma &= \frac{g}{V} (n-\cos \gamma),\quad &&\dot V &&= \frac{T-D}{m} - g \sin \gamma,\label{eq:sysDE2}
\end{alignat}
the initial conditions $x(0) = x_0$ and $z(0) = z_0$, the bound constraints
\begin{gather}
x(t+\delta t) > x(t),\quad \foralla \delta t > 0,\label{eq:Sx1}\\
V > 0,\quad 0 \le T \le T_{\max},\\
\left|\alpha\right| \le \frac{\pi}{18},\label{eq:AOTBC1}
\end{gather}
\end{subequations}
and the periodic constraints $z, \gamma, V \in \MBT_{T_f}$, where $\sigma$ is the Thrust Specific Fuel Consumption (TSFC), which gives the fuel efficiency of an engine design with respect to thrust output, $T_f \in \MBRP$ is the free final time, $g$ is the gravitational acceleration, $n$ is the load factor, $D$ is the drag, $T_{\max}$ is the thrust maximum value, and $\{x_0, z_0\} \subset \MBR$. The load factor and the drag are given by
\[n = \frac{L}{m g},\quad D = \frac{1}{2} \rho S \left[C_{D_0}+\frac{(C_{L0} + C_{L \alpha} \alpha)^2}{\pi\rme AR}\right] V^2,\]
where $L$ is the lift, $m$ is the UAV mass, $\rho$ is the air density, $S$ is the UAV wing planform area, $C_{D_0}$ is the constant parasitic drag, $C_{L0}$ and $C_{L \alpha}$ are the lift aerodynamic coefficients\footnote{$C_{L0}$ is the coefﬁcient of Lift at zero $\alpha$ and $C_{L \alpha}$ is the coefficient of $\alpha$ induced lift.}, $\rme$ is the Oswald efficiency factor, and $AR$ is the wing aspect ratio. The lift is further given by
\begin{equation}
L = \frac{1}{2} \rho S (C_{L0} + C_{L \alpha} \alpha) V^2.
\end{equation}
We refer to the OCP described in this section by OCP 1. Notice that Ineq. \eqref{eq:Sx1} implies $\dot x > 0\,\forall t \in \FOmega_{T_f}$. Moreover, the bound constraint \eqref{eq:AOTBC1} ensures a small attack angle to reduce the induced drag and maintain the validity of the UAV model. 

\section{FPS Interpolation and Quadrature}
\label{sec:FIPSQ1}
Let $\C{T} \in \MBRP, N \in \MBZeP, t_{j} \in \MBS_N^{\C{T}}\,\forall j \in \MBJ_N$, and $f_j = f(t_{j})\,\foralls f \in \MBT_{\C{T}}$. Through the DFT pair, one can write the $N/2$-degree, $\C{T}$-periodic Fourier interpolant, ${I_N}f$, of $f$ on the equispaced points grid $t_{0:{N-1}}$ as
\begin{equation}\label{eq:FI1nn1}
{I_N}f(t) = \sumd\sum\limits_{\left| k \right| \le N/2} {{\tilde f_k} {e^{i{\omega _k}t}}} = \sum\limits_{j = 0}^{N - 1} {{f_j}{\C{F}_j}(t)}, 
\end{equation}
where the primed sigma denotes a summation in which the last term is omitted, ${\tilde f_k}$ is the discrete Fourier interpolation coefficient given by
\[{\tilde f_k} = \frac{1}{N}\sum\limits_{j = 0}^{N - 1} {{f_j}{e^{ - i \omega_k {t_{j}}}}}\quad \forall k \in \MBK_N,\]
and ${\C{F}_j}(t)$ is the $N/2$-degree, $\C{T}$-periodic trigonometric Lagrange interpolating polynomial given by
\[{\C{F}_j}(t) = \frac{1}{N}\sumd\sum\limits_{\left| k \right| \le N/2} {e^{i{\omega _k}(t - {t_{j}})}} = {\left[ {\frac{1}{N}\sin \left( {\frac{{\pi N}}{\C{T}}\left( {t - {t_{j}}} \right)} \right)\cot \left( {\frac{\pi }{\C{T}}\left( {t - {t_{j}}} \right)} \right)} \right]_{t \ne {t_{j}}}},\]
$\forall j \in \MBJ_N$. One can approximate the definite integrals of $f$ over the successive intervals $\FOmega_{t_{0:N-1}}$ using Fourier PS quadrature as follows:
\begin{equation}\label{eq:AppRedFIM1}
{{\C{I}_{t_l}^{(t)}f}} \approx \sum\limits_{j = 0}^{N - 1} {{\theta^{\C{T}} _{l,j}}{f_j}},\quad \forall l \in \MBJ_N,
\end{equation}
where 
\begin{equation}\label{eq:RedFIM1}
{\theta^{\C{T}} _{l,j}} = \frac{1}{N}\left[ {{t_{l}} + \frac{{\C{T}i}}{{2\pi }} \sumd\sum\limits_{\scriptstyle \left| k \right| \le N/2\atop
\scriptstyle k \ne 0} {\frac{1}{k}{e^{ - i{\omega _k}{t_{j}}}}\left( {1 - {e^{i{\omega _k}{t_{l}}}}} \right)}} \right],\quad \forall l,j \in \MBJ_N,
\end{equation}
are the entries of the first-order square FIM, $\FThe^{\C{T}}$, of size $N$ and associated with the interval $\FOmega_{\C{T}}$. We can also describe the system of quadratures \eqref{eq:AppRedFIM1} in matrix notation as
\[{{\C{I}_{{{\bmt_N}}}^{(t)}f}} = \FThe^{\C{T}} f_{0:N - 1}.\]
In the special case when $\C{T} = 2 \pi$, Eq. \eqref{eq:RedFIM1} simply reduces to
\begin{equation}\label{eq:RedFIM12}
{\theta^{2\pi} _{l,j}} = \frac{1}{N}\left[ {{t_{l}} + i \sumd\sum\limits_{\scriptstyle \left| k \right| \le N/2\atop
\scriptstyle k \ne 0} {\frac{1}{k}{e^{ - i k {t_{j}}}}\left( {1 - {e^{i k {t_{l}}}}} \right)}} \right],\quad \forall l,j \in \MBJ_N.
\end{equation}
Eq. \eqref{eq:RedFIM12} saves a considerable amount of computational effort compared with Eq. \eqref{eq:RedFIM1} for large mesh grids. Besides, $\FThetpi$ is $\C{T}$-invariant, so we can precompute and store it first and then invoke it later quickly once we run the code. Therefore, $\FThetpi$ is practically the optimal FIM among the spectrum of $\C{T}$ values in terms of computational complexity, numerical stability, and speed. In addition, using the change of variables
\begin{equation}
t = \mu \tau:\quad \mu = \frac{\C{T}}{2\pi},
\end{equation}
we find that
\begin{gather}
\FThe^{\C{T}}_l f_{0:N - 1} = \C{I}_{t_l}^{(t)}{I_Nf} = \mu \C{I}_{\omega_{t_l}}^{(\tau)}{I_Nf\cancbra{\mu \tau}} = \mu \C{I}_{\tau_l}^{(\tau)}{I_Nf\cancbra{\mu \tau}}\\
= \mu \FThetpi_l I_Nf\left(\mu \bmtau_N\right) = \mu \FThetpi_l f_{0:N - 1}\quad \forall l \in \MBJ_N,
\end{gather}
where $t_l = \mu \tau_l\,\forall l \in \MBJ_N$ and $\FThe^{r}_l$ denotes the $l$th row of $\FThe^{r}\,\forall r > 0$. Hence, 
\begin{equation}\label{eq:Rel1}
\FThe^{\C{T}} = \mu \FThetpi.
\end{equation}
Although Eq. \eqref{eq:Rel1} allows us to directly generate $\FThe^{\C{T}}$ from $\FThetpi$ by premultiplying the latter by the scaling factor $\mu$, to approximate $\C{I}_{t_l}^{(t)}f$ in practice, it is more computationally effective to compute $\FThetpi_l f_{0:N - 1}$ first then multiply the result by $\mu$ rather than generating $\FThe^{\C{T}}$ first using Eq. \eqref{eq:Rel1} and then compute $\FThe^{\C{T}}_l f_{0:N - 1}$, since the former procedure requires $2 N^2$ flops, while the latter entails $N (3N-1)$ flops. We refer to $\FThetpi$ by the ``basic/principle'' FIM, because of its simplest form, the ``generating'' FIM due to Formula \eqref{eq:Rel1}, or the ``natural'' FIM, since its merits place it as the natural choice among all possible FIMs. For further information about FPS interpolation and quadratures, the reader may consult the recent works \cite{elgindy2019high,Elgindy2023a,Elgindy2023b}.

\section{The FIPS-ED Method}
\label{sec:FPIMIRF1}
Let $\alpha' = C_{L0} + C_{L \alpha} \alpha, c_1 = \sigma/(2 \pi), c_2 = \frac{1}{2} \rho S, c_3 = c_2/m$, and $c_4 = \pi \rme AR$. In order to use the natural FIM, we normalize the time horizon using the change of variables
\begin{equation}
t = \mu \tau: \mu = \frac{T_f}{2 \pi},
\end{equation}
and integrate the differential equations system \eqref{eq:sysDE1} and \eqref{eq:sysDE2} to obtain the following normalization of OCP 1 in integral form:
\begin{subequations}
\begin{equation}\label{eq:eq12}
\mathop {\min }\limits_{\tilalp, \tilT} J = c_1 \C{I}_{2\pi}^{(\tau)} {\tilT}
\end{equation}
subject to the 2D integral point-mass UAV dynamic model 
\begin{gather}
\tilx = \mu \C{I}_{\tau}^{(\tau)} {\left(\tilV \cos \tilgam\right)} + x_0,\quad \tilz = \mu \C{I}_{\tau}^{(\tau)} {\left(\tilV \sin \tilgam\right)} + z_0,\label{eq:sysDE12}\\
\tilgam = \mu \C{I}_{\tau}^{(\tau)} {\frac{c_3 \tilalp' \tilV^2-g \cos \tilgam}{\tilV}} + \tilgam(0),\label{eq:sysDE22}\\
\tilV = \mu \C{I}_{\tau}^{(\tau)}{\left(\frac{\tilT-\tilD}{m} - g \sin \tilgam\right)} + \tilV(0),\label{eq:sysDE23}
\end{gather}
the bound constraints
\begin{gather}
\cos \tilgam > 0,\quad \tilV > 0,\quad \left|\tilalp\right| \le \frac{\pi}{18},\quad 0 \le \tilT \le T_{\max},\quad T_f > 0,
\end{gather}
\end{subequations}
and the periodic constraints $\tilz, \tilgam, \tilV \in \MBT_{2\pi}$, where each shifted variable $\tilxi(\tau) = \xi(\mu \tau)\,\forall \xi$ $\in \digamma = \{x,z,\gamma,V,\alpha,\alpha',T,D\}$. We refer to this normalized integral-form OCP by OCP 2. Notice that all integrand functions of OCP 2 are $2\pi$-periodic, so their definite integrals can be effectively treated by the natural FIM. To reduce the computational cost of evaluating the performance index through numerical optimization procedures, we consider OCP 2 with the scaled performance index 
\begin{equation}
\barJ = \frac{N}{\sigma} J,\quad \foralls N \in \MBZeP,
\end{equation}
in lieu of $J$. This scaled version of OCP 2, denoted by OCP 3, is equivalent to OCP 2 in the sense that an optimal solution to OCP 2 is also an optimal solution to OCP 3. Fourier collocation of OCP 3 at the mesh points set $\MBS_N^{2\pi}$ converts the problem into the following constrained NLP:
\begin{mini}
   {\tilbmalp_N, \tilbmT_N}{\barJ_N \approx \bmone_N^{\top} \tilbmT_N}{}{}
   {\label{prob:Opt1}}{}
   \addConstraint{\tilbmx_N}{\approx \mu \FThetpi \left(\tilbmV_N \odot \cos \tilbmgam_N\right) + x_0 \bmone_N,}{}
   \addConstraint{\tilbmz_N}{\approx \mu \FThetpi \left(\tilbmV_N \odot \sin\tilbmgam_N\right) + z_0 \bmone_N,}{}
   \addConstraint{\tilbmgam_N}{\approx \mu \FThetpi \left[\left(c_3 \tilbmalp'_N \odot \left(\tilbmV_N\right)_2 - g \cos\tilbmgam_N\right) \oslash \tilbmV_N\right]}{}
   \addConstraint{}{+ \tilbmgam(0) \bmone_N,}{}
   \addConstraint{\tilbmV_N}{\approx \mu \FThetpi \left[\frac{1}{m} \left(\tilbmT_N-\tilbmD_N\right) - g \sin\tilbmgam_N\right] + \tilbmV(0) \bmone_N,}{}
   \addConstraint{\cos \tilbmgam_N }{> \bmzer_N,\quad \tilbmV_N > \bmzer_N,\quad \left|\tilbmalp_N\right| \le \frac{\pi}{18} \bmone_N,}{}
   \addConstraint{\bmzer_N}{ \le \tilbmT_N \le T_{\max} \bmone_N,\quad T_f > 0,}{}
 \end{mini} 
 where 
 \begin{align}
 \tilbmalp'_N &= C_{L0} \bmone_N + C_{L \alpha} \tilbmalp_N,\\
 \tilbmD_N &= c_2 \left[C_{D_0} \bmone_N +\frac{1}{c_4} \left(\tilbmalp'_N\right)_2\right] \odot \left(\tilbmV_N\right)_2.
 \end{align}
To solve the above problem for the collocation set of shifted state and control vectors, $\left\{\tilbmxi_N: \xi \in \digamma\right\}$, using standard numerical optimization methods, we can further rewrite the problem in terms of the $(6N+1)$-dimensional solution vector $\bmX_{6N+1} = [\tilbmx_N; \tilbmz_N; \tilbmgam_N$; $\tilbmV_N; \tilbmalp_N; \tilbmT_N; T_f]$ as follows:
\begin{mini}
   {\bmX_{4N:6N-1}}{\bmone_N^{\top} \bmX_{5N:6N-1}}{}{}
   {\label{prob:Opt2}}{}
   \addConstraint{\bmX_{0:N-1}}{\approx \mu \FThetpi \left(\bmX_{3N:4N-1} \odot \cos \bmX_{2N:3N-1}\right)}{}
   \addConstraint{}{+ x_0 \bmone_N,}  
   \addConstraint{\bmX_{N:2N-1}}{\approx \mu \FThetpi \left(\bmX_{3N:4N-1} \odot \sin \bmX_{2N:3N-1}\right)}{}
   \addConstraint{}{+ z_0 \bmone_N,}
   \addConstraint{\bmX_{2N:3N-1}}{\approx \mu \FThetpi \left[\left(c_3 \tilbmalp'_N \odot \left(\bmX_{3N:4N-1}\right)_2\right.\right.}{}
\addConstraint{}{\left.\left. - g \cos\bmX_{2N:3N-1}\right) \oslash \bmX_{3N:4N-1}\right]}
\addConstraint{}{+ \bmX_{2N} \bmone_N,}
   \addConstraint{\bmX_{3N:4N-1}}{\approx \mu \FThetpi \left[\frac{1}{m} \left(\bmX_{5N:6N-1}-\tilbmD_N\right)\right.}{}
\addConstraint{}{\left. - g \sin\bmX_{2N:3N-1}\right] + \bmX_{3N} \bmone_N,}      
   \addConstraint{\cos \bmX_{2N:3N-1} }{> \bmzer_N,\quad \bmX_{3N:4N-1} > \bmzer_N,}{}
   \addConstraint{\left|\bmX_{4N:5N-1}\right|}{\le \frac{\pi}{18} \bmone_N,}{}
   \addConstraint{\bmzer_N \le \bmX_{5N:6N-1}}{\le T_{\max} \bmone_N,\quad \bmX_{6N} > 0}
    \end{mini} 
 where 
 \begin{align}
 \mu &= \frac{1}{2\pi} \bmX_{6N},\quad \tilbmalp'_N = C_{L0} \bmone_N + C_{L \alpha} \bmX_{4N:5N-1},\\
 \tilbmD_N &= c_2 \left[C_{D_0} \bmone_N +\frac{1}{c_4} \left(\tilbmalp'_N\right)_2\right] \odot \left(\bmX_{3N:4N-1}\right)_2.
 \end{align}
NLP \eqref{prob:Opt2} is then solved recursively for the predicted solution vector $\bmX_{6N+1}$ starting with some initial mesh size, $\Nin \in \MBZP$, and an incremental increase, $\Ninc \in \MBZP$, until the lower bound on the size of the period, $T_f$, falls below a certain user tolerance $\varepsilon$. After obtaining the solution vectors set $\left\{\tilbmxi_N: \xi \in \digamma\right\}$, we shift it back onto the real time domain $\FOmega_{\C{T}}$ to obtain the desired collocation set $\left\{\bmxi_N: \xi \in \digamma\right\}$, and the state and control variables can be readily recovered at any time $t \in \FOmega_{\C{T}}$ by using Eq. \eqref{eq:FI1nn1}.

In the following section, we introduce an upgraded version of the edge-detection method of \citet{Elgindy2023a} to accurately reconstruct the periodic bang-bang thrust from the FPS data. This is a crucial step, because the Fourier interpolant of the thrust obtained by solving NLP \eqref{prob:Opt2} exhibits spurious oscillations near the jump discontinuities for increasing collocation mesh size that slow down its convergence, similar to the Gibbs phenomenon of Fourier series. The discontinuity points of the thrust are not known a priori; therefore, a key step in reconstructing the thrust is the accurate approximations of the jump discontinuity points locations. 

\subsection{The FPSED Method}
\label{subsec:TFEDM1}
We propose an edge-detection method motivated by the fact that the location of a jump discontinuity point $\xi$ is gradually squeezed between the locations of the sharp spikes of the Fourier interpolant graph near $\xi$ and eventually (almost) falls at the midpoint between the two abscissas whose ordinates are the peak and the bottom out of the two jagged lines enclosing $\xi$, as $N \to \infty$. The remainder of this section provides a detailed description of the method.

Let $I_Nf$ be the Fourier interpolant of a function $f \in \C{B}_{\C{T}}$ with any number of jump discontinuity points in $\F{\Omega}_{\C{T}}$. Let also $\MBS_M^{\C{T}} = \{y_{0:M-1}\}\,\foralll M \in \MBZ^+, d_{\max} = \indmax I_Nf(\bm{y}_M)$, and $d_{\min} = \indmin I_Nf(\bm{y}_M)$. To improve the accuracy of the Fourier interpolant extreme points, we extremize $I_Nf$ on the relatively small uncertainty intervals $[{y_{{d_{\min }} + \delta_{d_{\min },0} - 1}},{y_{{d_{\min }} - \delta_{d_{\min },M} + 1}}]$ and $[{y_{{d_{\max }} + \delta_{d_{\max },0} - 1}}$, ${y_{{d_{\max }} - \delta_{d_{\max },M} + 1}}]$ using the rapid Chebyshev PS line search method (CPSLSM) \cite{Elgindy2018optimization}. Let $I_N^{\max}f$ and $I_N^{\min}f$ be the approximate maximum and minimum values of $I_Nf$, respectively. Let also $\epsilon = \tilde \epsilon\;(I_N^{\max}f - I_N^{\min}f)\,\foralls \tilde \epsilon \in (0, 0.01]$, and consider the separation line $y = I_N^{\text{ave}}f = \frac{1}{2} \left(I_N^{\max}f + I_N^{\min}f\right)$ and the $\epsilon$-discontinuity feasible zone $\C{Z}_{\epsilon}^{\text{disc}}$, as described in \cite{Elgindy2023a}. We determine the first set of approximate discontinuity (AD) points, $\Xi = \left\{\tilxi_{1:L_1}\right\}\,\foralls L_1 \in \MBZP$, such that $\foralle k \in \MBN_{L_1}, \tilxi_k = y_j: I_Nf(y_j) \in \C{Z}_{\epsilon}^{\text{disc}}\,\foralls j \in J^1_{1:L_1}$, where $J^1_{1:L_1}$ is an index vector of integers such that $\left\{J^1_{1:L_1}\right\} \subset \MBJ_M^+$. To determine the locations of the remaining discontinuities, if exist, define the discrete auxiliary function, $I_N^{\text{aux}}f$, by
\begin{equation}\label{eq:INaux1}
I_N^{{\text{aux}}}f({y_{l}}) = \left\{ \begin{array}{l}
I_N^{\min }f,\quad \text{if }{I_N}f({y_{l}}) < I_N^{{\text{ave}}}f,\\
I_N^{\max }f,\quad \text{if }{I_N}f({y_{l}}) > I_N^{{\text{ave}}}f,
\end{array} \right.\quad \forall l \in \MBJ_{M}^+,
\end{equation}
and determine the index vector of candidate discontinuity points $J^2_{1:L_2}\,\foralls L_2 \in \MBZeP$: 
\begin{equation}\label{eq:INaux1}
I_N^{{\text{aux}}}f(y_{J^2_l}) - I_N^{{\text{aux}}}f(y_{J^2_l+1}) \ne 0\quad \forall l \in \MBN_{L_2}.
\end{equation}
An initial data screening is recommended at this stage to avoid having clone, contiguous, or sufficiently close candidate discontinuity points from the obtained two index vectors in a relatively small vicinity of $\FOmega_{\C{T}}\,\foralll M$ value. This can be achieved by eliminating all indices in $\left\{J^2_{1:L_2}\right\}$ within a sufficiently close distance from $\left\{J^1_{1:L_1}\right\}$; i.e., we select a relatively small $r_1 \in \MBZP$ and set up the filtered index vector $J^3_{1:L_3}$ such that
\[\left\{J^3_{1:L_3}\right\} = \left\{j \in J^2_{1:L_2}: \dist\left(\left\{J^1_{1:L_1}\right\},\left\{J^2_{1:L_2}\right\}\right) > r_1\right\}\quad \foralls L_3 \in \MBZP.\]
The remaining discontinuity points $\xi_{L_1+1:L_1+L_3}$ either exist in the small uncertainty intervals $\left. {}_{y_{J_i}}\FOmega_{y_{J_i+1}}\right|_{i=1:L_3}$ or occur in close proximity of their boundaries; we estimate their values by the midpoints of these intervals. In particular, if we denote the boundary points $y_{J^3_{1:L_3}}$ and $\left. y_{J^3_i+1} \right|_{i=1:L_3}$ in ascending order by $b_{1:2 L_3}$, respectively, then we can readily calculate the midpoints of the uncertainty intervals $\left. {}_{b_{2 i-1}}\FOmega_{b_{2 i}}\right|_{i=1:L_3}$ and arrive at the following estimates:
\[\tilxi_{L_1+l} = \frac{1}{2} (b_{2 l-1}+b_{2 l})\quad \forall l \in \MBN_{L_3}.\]
This gives the updated set of AD points $\Xi := \Xi \cup \left\{\tilxi_{L_1+l:L_1+L_3}\right\} = \{\hxi_{1:L_2}\}$, where $\hxi_{1:L_2}$ are the AD points in ascending order. As a further correction step in practice, we set 
\begin{equation}\label{eq:hxi1}
\hxi_{L_2} := \C{T},
\end{equation}
if $\hxi_{L_2}$ is within a sufficiently small distance from $t = \C{T}$, which can be achieved by choosing a relatively small $r_2 \in \MBZP$, and then applying Eq. \eqref{eq:hxi1} when $J_{L_2} =  M - r_2\,\foralll M$; this step completes the search procedure. 

After completing the search procedure for locating the AD points, we need to damp the ripples of the Fourier interpolant graph adjacent to the extreme values lines $y = f_{\max}$ and $y = f_{\min}$, where $f_{\max}$ and $f_{\min}$ are the maximum and minimum values of $f$, respectively. To this end, we use the median to assess the central tendency of Fourier interpolant values. In particular, let $\MB{FI}^{\text{u}}$ and $\MB{FI}^{\text{d}}$ denote the Fourier interpolant values sets above and below the separation line, and calculate the medians of these two data sets, denoted by ${I_{N,M}}f_+^{\text{med}}$ and ${I_{N,M}}f_-^{\text{med}}$ in respective order. Then we can define the corrected Fourier interpolant, denoted by $f^c_{N,M}$, by
\begin{equation}\label{eq:RecSt1}
f^c_{N,M}(t) = \left\{ \begin{array}{l}
{I_{N,M}}f_+^{\text{med}},\quad \text{if }t \in \FOmega_{\tilxi_1} \cup \left.{}_{\tilxi_{2l}}\FOmega_{\tilxi_{2l+1}}\right|_{l=1:\frac{L_2}{2}-1} \cup\;{}_{\tilxi_{L_2}}\FOmega_{\C{T}},\\
{I_{N,M}}f_-^{\text{med}},\quad \text{otherwise},
\end{array} \right.
\end{equation}
if ${I_N}f(0) \in \MB{FI}^{\text{u}}$, or by
\begin{equation}\label{eq:RecSt2}
f^c_{N,M}(t) = \left\{ \begin{array}{l}
{I_{N,M}}f_-^{\text{med}},\quad \text{if }t \in \FOmega_{\tilxi_1} \cup \left.{}_{\tilxi_{2l}}\FOmega_{\tilxi_{2l+1}}\right|_{l=1:\frac{L_2}{2}-1} \cup\;{}_{\tilxi_{L_2}}\FOmega_{\C{T}},\\
{I_{N,M}}f_+^{\text{med}},\quad \text{otherwise},
\end{array} \right.
\end{equation}
otherwise. If the extreme values of $f$ are known a priori, then the medians computations in the final reconstruction step are dispensed, and ${I_{N,M}}f_+^{\text{med}}$ and ${I_{N,M}}f_-^{\text{med}}$ in Eqs. \eqref{eq:RecSt1} and \eqref{eq:RecSt2} are replaced by $f_{\max}$ and $f_{\min}$, respectively.

\begin{figure}[H]
\centering
\hspace{-5mm}\includegraphics[scale=0.6]{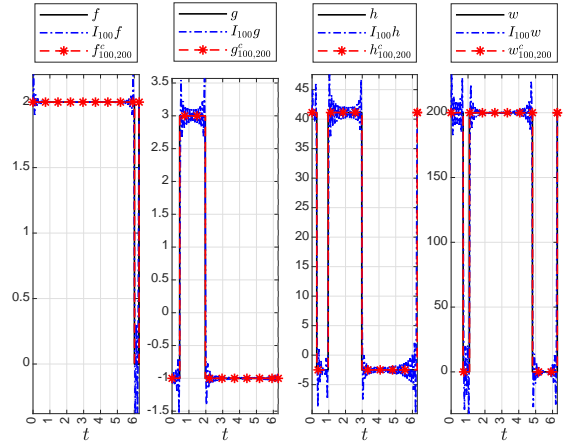}
\caption{Plots of $f, g, h$, and $w$ on $\FOmega_{2\pi}$ and their corresponding Fourier interpolants and corrections obtained by the FPSED method using $(N,M,r_{1:2}) = (100,200,1,2)$.}
\label{fig:Fig1}
\end{figure}

Figure \ref{fig:Fig1} shows a demonstration of how the FPSED method successfully reconstructs the following four periodic bang-bang test functions on $\FOmega_{2\pi}$ from their FPS data:\\
\scalebox{0.9}{\parbox{\linewidth}{%
\begin{gather}
f(t) = \left\{ \begin{array}{l}
2,\quad 0 \le t < 6.01,\\
0,\quad 6.01 \le t < 2\pi,\\
2,\quad t = 2\pi,
\end{array} \right.\quad g(t) = \left\{ \begin{array}{l}
 - 1,\quad 0 \le t < 0.45,\\
3,\quad 0.45 \le t \le 1.97,\\
 - 1,\quad 1.97 \le t \le 2\pi ,
\end{array} \right.\quad \\
h(t) = \left\{ \begin{array}{l}
41.12,\quad 0 \le t < 0.28,\\
 - 2.5,\quad 0.28 \le t \le 0.96,\\
41.12,\quad 0.96 \le t \le 2.98,\\
 - 2.5,\quad 2.98 \le t < 2\pi,\\
41.12,\quad t = 2\pi, 
\end{array} \right.\; w(t) = \left\{ \begin{array}{l}
200,\quad 0 \le t < 0.71,\\
0,\quad 0.71 \le t \le 1.08,\\
200,\quad 1.08 \le t \le 4.81,\\
0,\quad 4.81 \le t < 2\pi ,\\
200,\quad t = 2\pi .
\end{array} \right.
\end{gather}
}}\\
The FPSED method was implemented for all test functions using $(N,M,r_{1:2}) = (100,200,1,2)$, and the maximum absolute errors in computing $\tilxi_{1:L_2}$ are reported in Table \ref{tab:1}.

\begin{table}[H]
\centering
\caption{The AD points of the test functions obtained using the FPSED method with $(N,M,r_{1:2}) = (100,200,1,2)$. The shown approximations are rounded to five significant digits.}
\resizebox{0.45\textwidth}{!}{  
\begin{tabular}{cccccc}  
\toprule
Test function & $\tilxi_1$ & $\tilxi_2$ & $\tilxi_3$ & $\tilxi_4$ & $MAE$\\
\midrule
$f$ & $6.0148$ & $6.2832$ & & & $4.8080e-03$\\
$g$ & $0.4578$ & $1.9734$ & & & $7.8200e-03$\\
$h$ & $0.2842$ & $0.9630$ & $2.9837$ & $6.2832$ & $4.1642e-03$\\
$w$ & $0.7104$ & $1.0893$ & $4.8150$ & $6.2832$ & $9.2959e-03$\\
\bottomrule
\end{tabular}
}
\label{tab:1}
\end{table}

\begin{rem}
Notice that $L_1 \le L_2$ in the proposed FPSED method, since the Fourier interpolant graph moves in virtually a vertical direction through each point whose abscissa is a discontinuity point, so that the Fourier interpolant values at the interpolation points that are not sufficiently close to a discontinuity point live outside the user-defined $\C{Z}_{\epsilon}^{\text{disc}}$, for relatively small $\epsilon$; thus, we expect to have at most a single mesh point near each discontinuity point, while $L_2$ is a positive even integer equal to the total number of discontinuities, because the periodic bang-bang function switches abruptly once at each discontinuity; consequently, Condition \eqref{eq:INaux1} must occur once near each discontinuity.
\end{rem}

\section{Errors and Convergence of Fourier Interpolation and Qaudrature} 
\label{sec:CRFSPF1}
For completeness sake, we state the following two theorems which can be found in \cite[Corollaries 5.1 \& 5.2]{Elgindy2023a,Elgindy2023b}.

\begin{thm}[Fourier interpolation and quadrature errors for analytic, $\C{T}$-periodic functions]\label{thm:2}
Suppose that $f \in {{\C{A}_{\C{T},\beta}}} \foralls \beta > 0$ is approximated by the $\C{T}$-periodic Fourier interpolant $I_Nf\,\foralls$ $N \in \MBZeP$, then $\left\| f - {{I_Nf}} \right\| = O\left( {{e^{ - {\omega _{N\beta/2}}}}} \right)$ and
\begin{subequations}
\begin{equation}\label{eq:FTS1AE1hi1}
\left| {\C{I}_{{{\bmt_N}}}^{(t)}f - \FThe^{\C{T}} f_{0:N - 1}} \right| = O\left( {{e^{ - {\omega _{N\beta/2}}}}\bmone_N} \right),\quad \text{as }N \to \infty.
\end{equation}
Moreover, if $f$ is $\beta$-analytic, then 
\begin{equation}\label{eq:Thm2hi2}
{\left\| {f - {I _N}f} \right\|} = 0\quad\text{and}\quad \left| {\C{I}_{{{\bmt_N}}}^{(t)}f - \FThe^{\C{T}} f_{0:N - 1}} \right| = \bmzer_N.
\end{equation}
\end{subequations}
\end{thm}

\begin{thm}[Fourier interpolation and qaudrature errors for nonsmooth and $\C{T}$-periodic functions]\label{thm:1}
Suppose that $f \in {{\C{H}_{\C{T}}^s}} \foralls$ $s \in \MBZ^+$ is approximated by the $\C{T}$-periodic Fourier interpolant $I_Nf\,\foralls$ $N \in \MBZeP$, then $\left\| f - {{I_Nf}} \right\| = O\left( {{N^{-s-1/2}}} \right)$ and 
\begin{equation}\label{eq:FTS1AE1hi126Apr2021}
\left| {\C{I}_{{{\bmt_N}}}^{(t)}f - \FThe^{\C{T}} f_{0:N - 1}} \right| = O\left( {{N^{-s-1/2}}} \bmone_N \right),\quad \text{as }N \to \infty.
\end{equation}
\end{thm}

Furthermore, Fourier interpolant $I_Nf\,\foralls N \in \MBZeP$ converges to a function $f \in {\C{H}_{\C{T}}^0}$ at the rate $O(1/N)$ for each continuity point of the function and oscillates at $O(1)$ on a small vicinity of a discontinuity point.

\section{Numerical Simulations}
\label{sec:NS1}
\begin{figure}[H]
\centering
\includegraphics[scale=0.25]{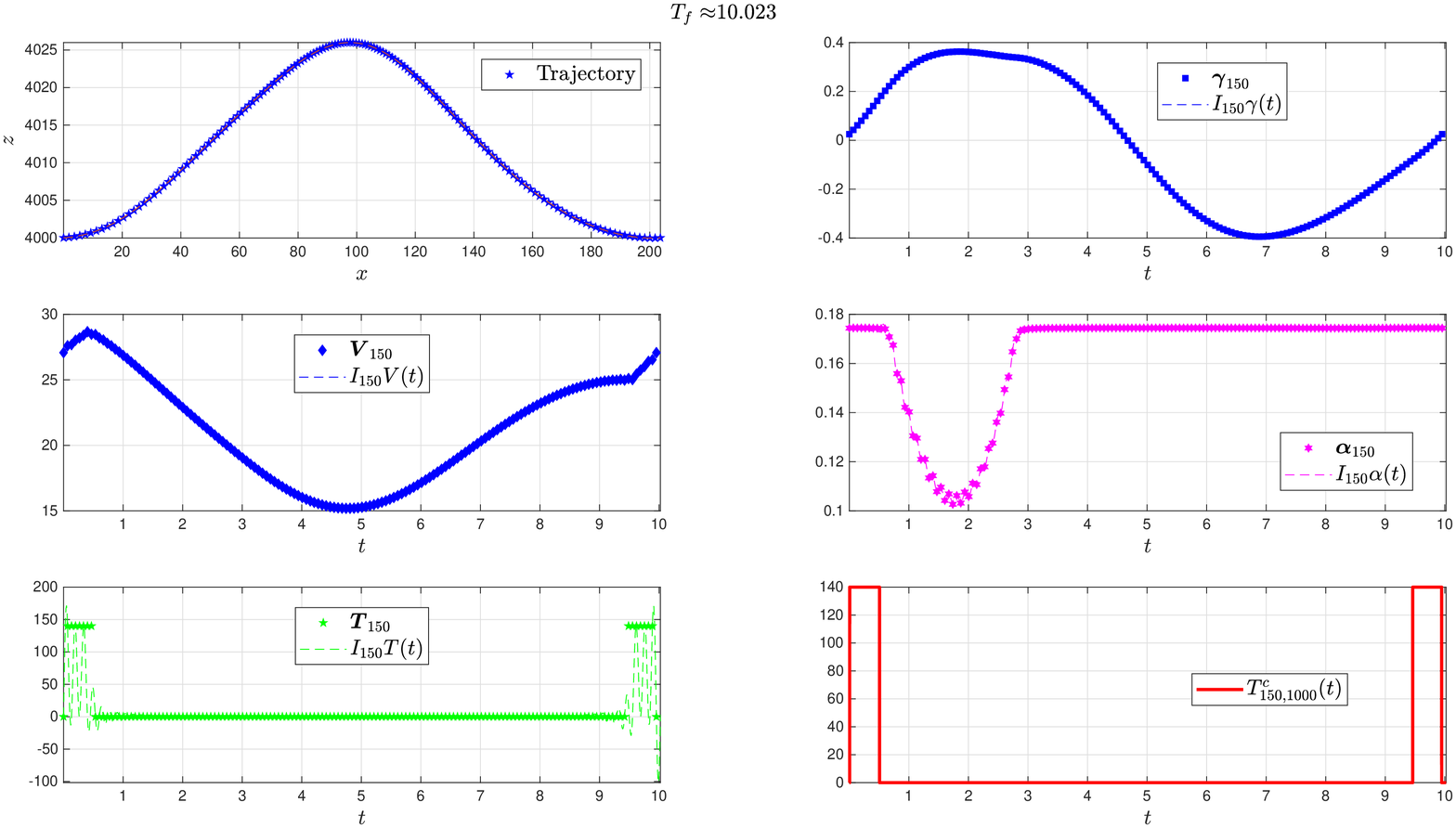}
\caption{The plots of the UAV flight trajectory (upper left), flight-path angle (upper right), speed (middle left), attack angle (middle right), predicted thrust (bottom left), and the corrected thrust (bottom right) profiles obtained using the FIPS-ED method with $N = 150$ and the parameter values listed in Table \ref{tab:3}. The calculated period is approximately $10.023\,s$.}
\label{fig:Fig2}
\end{figure}

\begin{figure}[H]
\centering
\includegraphics[scale=0.25]{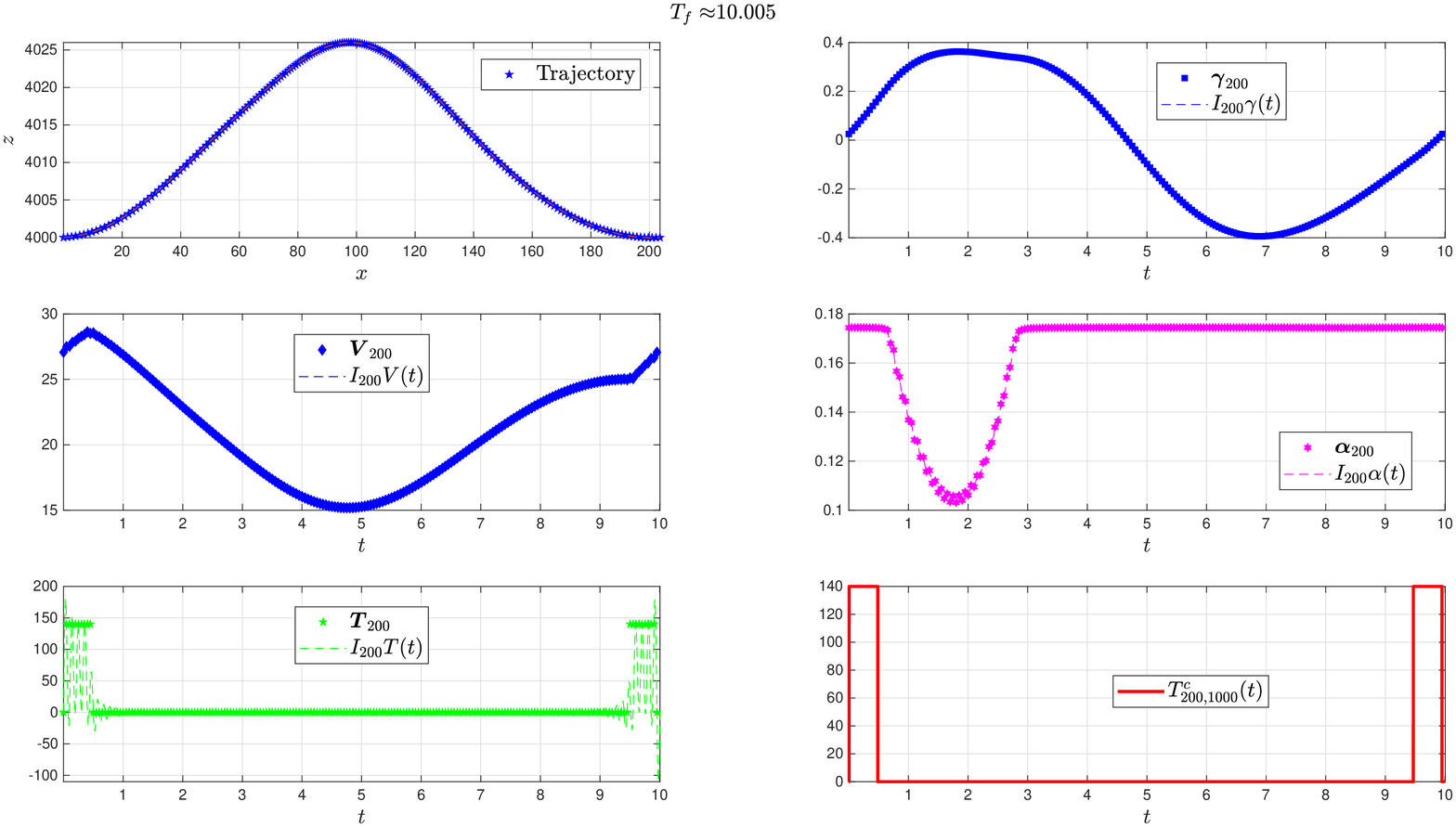}
\caption{The plots of the UAV flight trajectory (upper left), flight-path angle (upper right), speed (middle left), attack angle (middle right), predicted thrust (bottom left), and the corrected thrust (bottom right) profiles obtained using the FIPS-ED method with $N = 200$ and the parameter values listed in Table \ref{tab:3}. The calculated period is approximately $10.005\,s$.}
\label{fig:Fig3}
\end{figure}

\begin{figure}[H]
\centering
\includegraphics[scale=0.25]{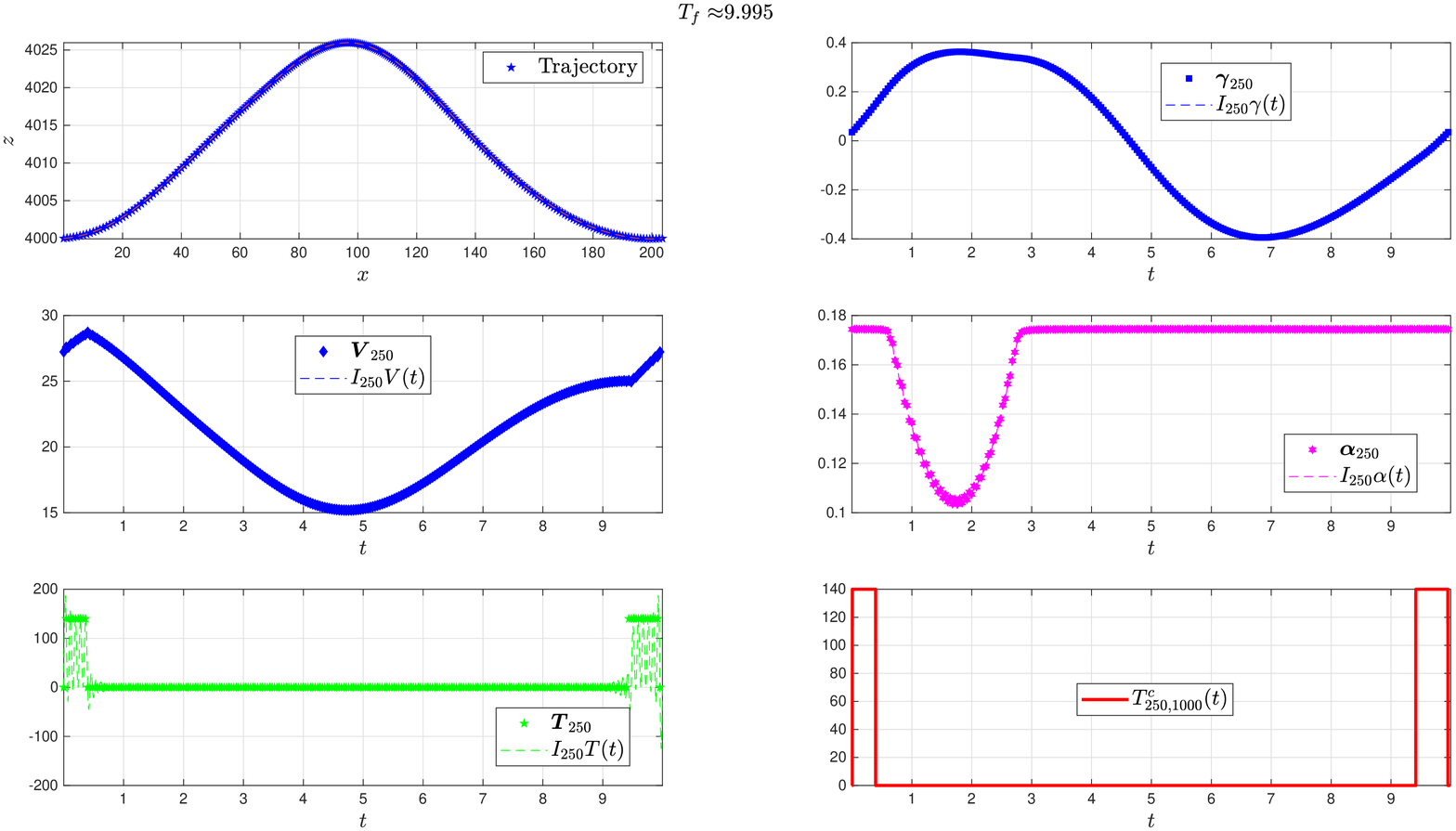}
\caption{The plots of the UAV flight trajectory (upper left), flight-path angle (upper right), speed (middle left), attack angle (middle right), predicted thrust (bottom left), and the corrected thrust (bottom right) profiles obtained using the FIPS-ED method with $N = 250$ and the parameter values listed in Table \ref{tab:3}. The calculated period is approximately $9.995\,s$.}
\label{fig:Fig4}
\end{figure}

\begin{figure}[H]
\centering
\includegraphics[scale=0.25]{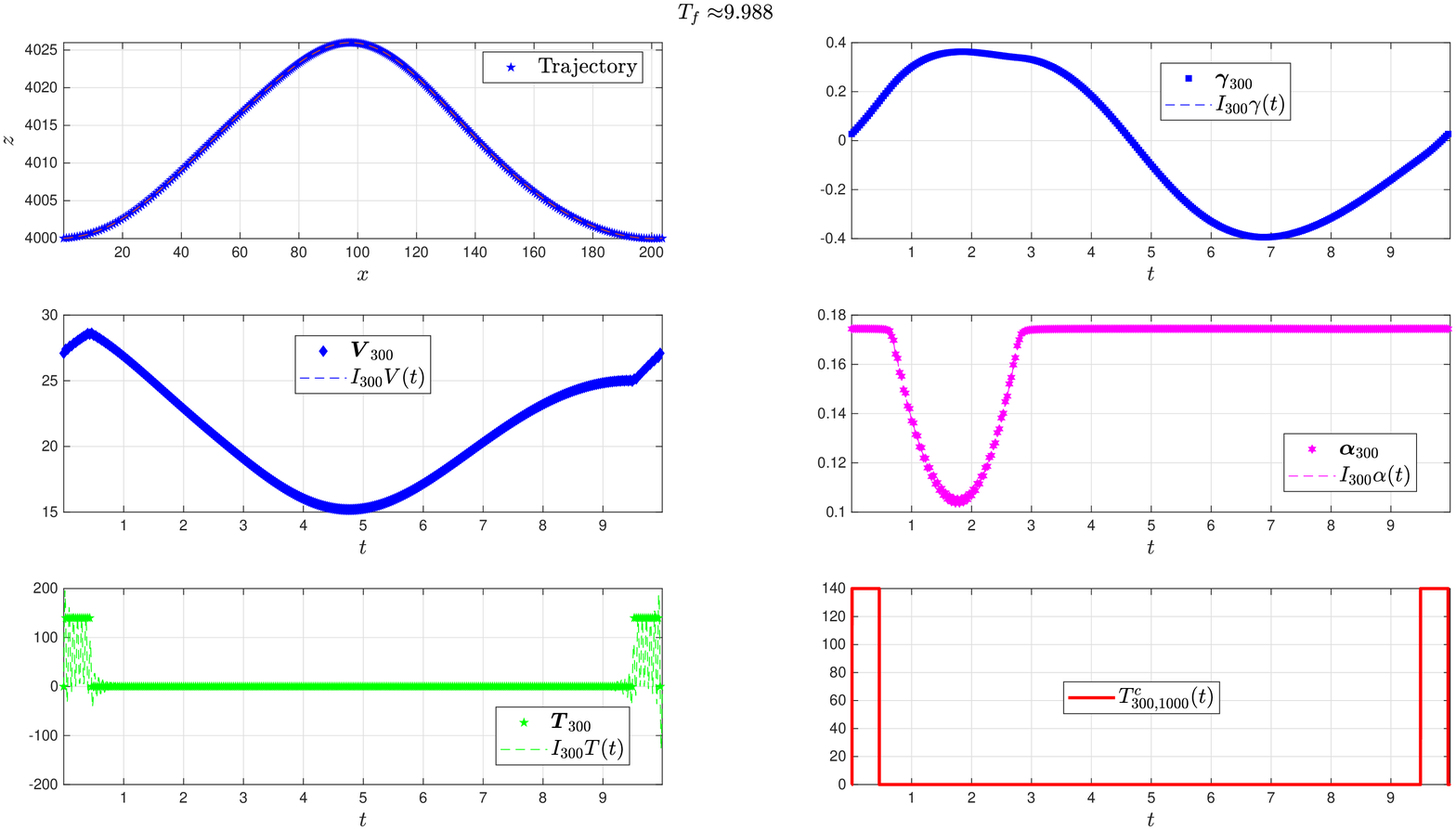}
\caption{The plots of the UAV flight trajectory (upper left), flight-path angle (upper right), speed (middle left), attack angle (middle right), predicted thrust (bottom left), and the corrected thrust (bottom right) profiles obtained using the FIPS-ED method with $N = 300$ and the parameter values listed in Table \ref{tab:3}. The calculated period is approximately $9.988\,s$.}
\label{fig:Fig5}
\end{figure}

This section presents some numerical experiments to demonstrate the performance of the proposed FIPS-ED method. The numerical experiments were carried out using MATLAB R2022b software installed on a personal laptop equipped with a 2.9 GHz AMD Ryzen 7 4800H CPU and 16 GB memory running on a 64-bit Windows 11 operating system. The FIPS-ED method was performed together with the MATLAB fmincon solver using the interior-point algorithm and the input parameter values listed in Table \ref{tab:3}. Starting with an initial mesh size of $150$ equally-spaced points with a mesh size incremental increase of $50$ points, the FIPS-ED method converges to an approximate optimal fuel consumption of $0.154\,kg/s$ when $N = 300$ with the approximate period $9.99\,s$. The optimal switching times of the thrust are approximately $\tilxi_{1:4} \approx 0.005\,s, 0.460\,s, 9.493\,s$, and $9.953\,s$, rounded to three decimal digits. The values of the initial flight-path angle and UAV speed are approximately $\gamma_0 \approx 0.027$ and $V_0 \approx 27.1 m/s$. The profiles of the state and control variables derived by the FIPS-ED method using $N = 150:50:300$ are shown in Figures \ref{fig:Fig2}-\ref{fig:Fig5}. Roughly speaking, over a period of about $10\,s$, the figures show that the UAV should swiftly produce a maximum thrust of $140\,N$ for nearly half a second to speed up the vehicle from about $27\,m/s$ at an altitude of $4,000\,m$ into a maximum speed of nearly $28.6\,m/s$ while climbing upward. The UAV propellers should produce zero thrust for nearly the next $9\,s$. During this period, the UAV speed initially declines gradually until it reaches a minimum of approximately $15.2\,m/s$ at a maximum altitude of about $4,026\,m$ and then steadily increase as the vehicle descents. In the last half of a second, the thrust should return to its maximum value, and the vehicle accelerates at a faster rate until it reaches its initial altitude. In summary, if we consider the entire flight endurance relative to the thrust history, then after approximately $9.5\,s$ of the first cycle, the thrust should practically be at its peak for nearly $1\,s$ and the propellers should produce zero thrust for the next $9\,s$, for each new cycle. During the entire flight, the magnitude of the flight-path angle is bounded by $0.4$. Moreover, if we denote the total number of cycles in the UAV endurance by $N_f$, then the attack angle remains nearly constant at a maximum value of approximately $0.174$, except for the approximate time intervals ${\left. {{}_{0.6 + n{T_f}}{\FOmega _{2.9 + n{T_f}}}} \right|_{n = 0:{N_f-1}}}$, where it glides down and up, arriving at a minimum value of about $0.1$ at $t \approx 1.8+n T_f\,\forall n = 0:N_f-1$.

\begin{figure}[H]
\centering
\includegraphics[scale=0.25]{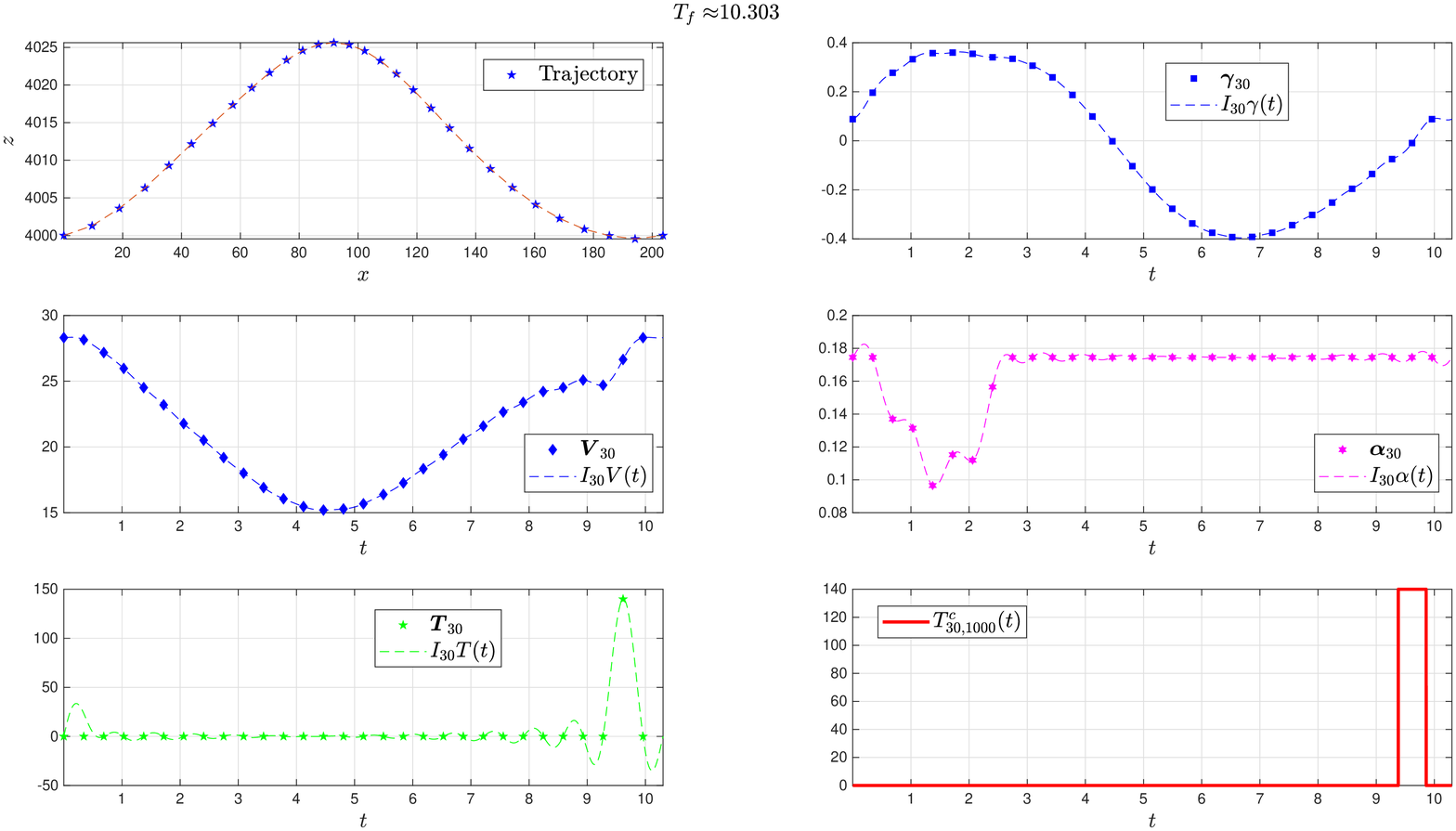}
\caption{The plots of the UAV flight trajectory (upper left), flight-path angle (upper right), speed (middle left), attack angle (middle right), predicted thrust (bottom left), and the corrected thrust (bottom right) profiles obtained using the FIPS-ED method with $N = 30$ and the parameter values listed in Table \ref{tab:3}. The calculated period is approximately $10.303\,s$.}
\label{fig:Fig6}
\end{figure}

\begin{figure}[H]
\centering
\includegraphics[scale=0.25]{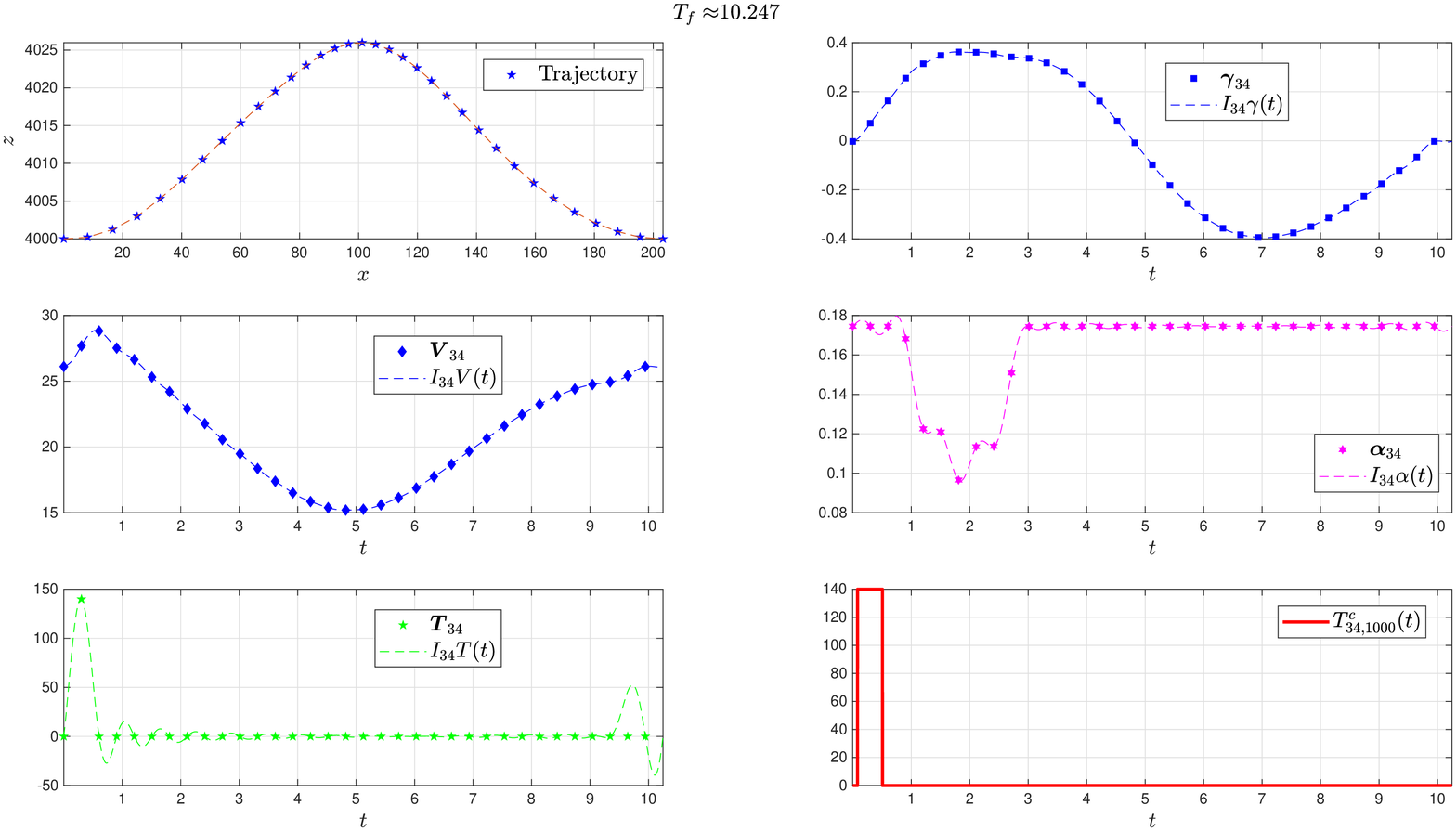}
\caption{The plots of the UAV flight trajectory (upper left), flight-path angle (upper right), speed (middle left), attack angle (middle right), predicted thrust (bottom left), and the corrected thrust (bottom right) profiles obtained using the FIPS-ED method with $N = 34$ and the parameter values listed in Table \ref{tab:3}. The calculated period is approximately $10.247\,s$.}
\label{fig:Fig7}
\end{figure}

It should be noted that the approximate optimal fuel consumption obtained by the FIPS-ED method is larger than those recorded in \cite{ogunbodede2019endurance} using the shooting method and the differential flatness based method, which were $0.07686\,kg/s$ and $0.07710\,kg/s$ with periods $9.96\,s$ and $10.08\,s$, respectively. However, there are a number of observations which largely support our results; we describe them in the following items:
\begin{itemize}
\item The basis function used in the differential flatness based method does not allow bang-bang profiles due to continuity constraints as reported in \cite{ogunbodede2019endurance}, so the bang-bang structure of the thrust profile is generally not preserved, and the method may skip some switching times of the optimal thrust policy.
\item A saturation function was wrapped around the values of the thrust control in \cite[Eq. (40)]{ogunbodede2019endurance}, which may affect the thrust structure in practice. In particular, the proposed formula 
\begin{equation}\label{eq:SF1}
T(t) := \Phi_{s_m}(t) = 70 + \frac{{140}}{\pi }{\tan ^{ - 1}}\left( { - {s_m}\left( {\sigma  + \frac{{{\lambda _V(t)}{T_f}}}{m}} \right) + \frac{\pi }{2}} \right),
\end{equation}
includes the usual parameters $\sigma$ and $m$, the final time $T_f$, the time-dependent Lagrange multiplier associated with the UAV speed $\lambda_V$, and a user-defined constant smoothing factor $s_m$ of large value, which varies depending on how sharp the bang-bang thrust is needed in the sense that $\mathop {\lim }\limits_{{s_m} \to \infty } {\Phi _{{s_m}}}(t) = T(t)\,\forall t \in \FOmega_{T_f}$. However, $\Phi _{{s_m}}$ is very ill-conditioned for $\lambda_V(0)$ values near the recorded numerical value $-0.01626$ in \cite{ogunbodede2019endurance} with increasing sensitivity as $s_m \to \infty$. Table \ref{tab:4} shows the values of $\Phi _{{s_m}}$ for $\lambda_V(0) = -0.01626$ and its perturbed value $\bar \lambda_V(0) = -0.01627$ for some increasing values of $s_m$; the values of $\sigma$ and $m$ are as quoted from Table \ref{tab:2}. Notice how the data error propagates when $\Phi _{{s_m}}$ is evaluated at $\bar \lambda_V(0)$ instead of the recorded $\lambda_V(0)$. For example, a small perturbation of $10^{-5}$ in $\lambda_V(0)$ when $s_m = 10^8$ causes a relative error in the calculation of $\Phi _{{s_m}}$ of approximately $198.5 \%$, rounded to four significant digits; therefore, the relative change in evaluating $\Phi _{{s_m}}$ is about $3,227$ times larger than the relative change in $\lambda_V(0)$(!) This shows that the thrust may not have been evaluated accurately near $t = 0$ using finite-digit arithmetic.
\item Table \ref{tab:4} shows that $\Phi _{{s_m}}$ is also sensitive with respect to $s_m$ values, dropping, for example, from $109$ into $19.3$ for $\lambda_V(0) = -0.01626$ when $s_m$ increases by one order of magnitude. Adaptivity should be exercised here with caution to determine the feasible range of $s_m$ values required for stability and to prevent possible divergence in practice.
\item Since $\lambda_V$ is periodic based on \cite[Eq. (27)]{ogunbodede2019endurance}, $\Phi _{{s_m}}$ must be periodic with respect to $\lambda_V$ when holding all other parameters fixed, thus, $T(0) = T(T_f)$ by Eq. \eqref{eq:SF1}, which does not conform to the thrust profile in \cite[Figure 4]{ogunbodede2019endurance}.
\item The work in \cite{ogunbodede2019endurance} aims at finding optimal periodic solutions to minimize the fuel consumption in UAVs to enhance endurance. However, the proposed definition \eqref{eq:SF1} of the thrust $T$ includes the nonperiodic inverse tangent function, which eliminates the periodicity property of $T$.
\item Fourier series and Fourier interpolants generally behave similar to each other for closely similar number of modes and meshes, except possibly at the discontinuity points. In particular, Fourier series converges to $\left(\xi^- + \xi^+\right)/2$ at a jump discontinuity point, while Fourier interpolant value falls within the open interval $\left(f(\xi^-), f (\xi^+)\right)$, except when the jump discontinuity point $\xi$ coincides with an interpolation node, where the Fourier interpolant matches the value of the discontinuous function according to the interpolation condition \cite{Elgindy2023a}, where $\xi^-$ and $\xi^+$ are points infinitesimally to the left and right of a discontinuity point $\xi$, respectively. With this in mind, and the fact that the differential flatness based method was performed in \cite{ogunbodede2019endurance} using only $40$ Fourier terms, while the FIPS-ED method was carried out using large collocation meshes to resolve the thrust bang-bang structure accurately, we expect the FIPS-ED method to probably overlook some switching times of the optimal bang-bang thrust policy for low mesh densities. Indeed, we observe in Figure \ref{fig:Fig6} that the FIPS-ED method attempts to resolve the abrupt switches near the end of the time horizon $\FOmega_{T_f}$ using $30$ mesh points, but fails to recognize the abrupt changes near the start. An opposite behavior is observed in Figure \ref{fig:Fig7}, where the method identifies approximate abrupt switches near $t = 0$ using $34$ mesh points, but again fails to resolve the shock near $t = T_f$. This shows that more mesh points in Fourier collocation have to be clustered near the boundaries of $\FOmega_{T_f}$ and equivalently more modes in Fourier series are needed to capture the thrust profile accurately.
\item The thrust profile looks almost the same for $N = 150:50:300$ indicating convergence of the FIPS-ED method to the optimal switching times within satisfactory accuracy.
\end{itemize}
The above arguments demonstrate that the smoothing technique in \cite{ogunbodede2019endurance} smears some of the discontinuities in the thrust, while the FIPS-ED method works well in capturing the abrupt changes in the thrust controller. The work in \cite{ogunbodede2019endurance} remains very important though, since it presents the first- and second-order conditions of optimality and provides an early numerical insight on the energy-optimal path planning problem governed by the 2D point-mass UAV dynamic model used to find the optimal periodic solutions to enhance the endurance of UAVs. 

\begin{table}[H]
 \centering
 \caption{$\Phi _{{s_m}}$ values for $\lambda_V(0) = -0.01626$ and $\bar \lambda_V(0) = -0.01627$ using $\sigma = 0.012, m = 13.5$, and some range of $s_m$ values. All calculations were rounded to one decimal digit.}
 \resizebox{0.35\textwidth}{!}{  
 \begin{tabular}{ccc}  
 \toprule
 $s_m$ & Recorded $\lambda_V(0)$ in \cite{ogunbodede2019endurance} & Perturbed $\lambda_V(0)$\\
 & $(-0.01626)$ & $(-0.01627)$ \\ 
 \midrule
 $1e5$ & $109.0$ & $118.7$ \\
 $1e6$ & $19.3$ & $131.6$ \\
 $1e7$ & $1.2$ & $138.8$ \\
 $1e8$ & $0.1$ & $139.9$ \\
 $1e9$ & $0.0$ & $140.0$ \\
 $1e10$ & $0.0$ & $140.0$ \\
 \bottomrule
 \end{tabular}
 }
 \label{tab:4}
 \end{table} 

\begin{table}[H]
\centering
\caption{Physical modeling parameters}
\resizebox{0.2\textwidth}{!}{  
\begin{tabular}{cc}  
\toprule
Parameters & Value\\
\midrule
$\sigma$ & $0.012\,kg/(s N)$ \\
$\rho$ & $1.2682\,kg/m^3$ \\
$S$ & $0.55\,m^2$ \\
$C_{D0}$ & $0.03$ \\
$C_{L0}$ & $0.28$ \\
$C_{L\alpha}$ & $3.45$ \\
$m$ & $13.5\,kg$ \\
$g$ & $9.81\,m/s^2$ \\
$\rme$ & $0.9$ \\
$AR$ & $15.2445$ \\
$T_{\max}$ & $140$ \\
$x_0$ & $0$\\
$z_0$ & $4000$\\
\bottomrule
\end{tabular}
}
\label{tab:2}
\end{table}

\begin{table}[H]
\centering
\caption{FIPS-ED method input parameters}
\resizebox{0.3\textwidth}{!}{  
\begin{tabular}{cc}  
\toprule
FIPS method parameters & Values\\
\midrule
$\Nin$ & $150$ \\
$\Ninc$ & $50$ \\
$M$ & $1000$ \\
$\varepsilon$ & $0.01$ \\
\midrule
FPSED method parameters & Values\\
$r_{1:2}$ & $1, 2$ \\
\midrule
fmincon parameters & Values\\
TolFun & $1e-12$ \\
TolX & $1e-12$ \\
\bottomrule
\end{tabular}
}
\label{tab:3}
\end{table}

\section{Conclusion}
\label{sec:Conc}
In this study, the energy-optimal path planning problem governed by the 2D point-mass UAV dynamic model is investigated numerically using Fourier collocation based on equally spaced meshes and Fourier quadrature induced by the accurate and very efficient natural FIM. The flight performance of the UAV was studied, including the propeller thrust characteristics, maneuverability, flight envelope, and flight attack angle, using input data and physical parameter values collected from \cite{ogunbodede2019endurance}. The convergence of the numerical tools necessary for the discretization of the OCP is addressed in detail. The FIPS-ED method readily converts the OCP into an NLP, which can be solved in the physical space for the state and control variables using a standard optimization software. Under certain parameter settings, our numerical study suggests that the UAV propellers should produce maximum thrust twice, rather than once as proposed in \cite{ogunbodede2019endurance}, for two short time periods during the flight operation to enhance UAV endurance while maintaining the validity of the UAV model. The proposed FPSED method proved to be a faithful companion to the applied FIPS method, which is able to accurately resolve the bang-bang profile of the approximate optimal thrust policy from the FPS data collected by grouping the FIPS method with MATLAB fmincon solver to maximize the UAV endurance. This provides the motivation for future studies on nonsmooth OCPs arising in various areas and applications using the developed techniques.  

\section*{Acknowledgment}
Special thanks to Dr. Oladapo Ogunbodede\footnote{Dr. Oladapo Ogunbodede is a control engineer at ThorDrive Inc, Cincinnati, OH, USA.} for valuable discussions and sharing the parameter values of $\sigma$ and $s_m$ used in this paper, which were not explicitly reported in \cite{ogunbodede2019endurance}.

\bibliographystyle{model1-num-names}
\bibliography{Bib}
\end{document}